\documentclass{amsart}

\newcommand{\dd}{\mathrm{d}}
\newcommand{\supp}{\mathop{\mathrm{supp}}}
\newcommand{\RR}{\mathbb{R}}

\newtheorem{theorem}{Theorem}
\newtheorem{lemma}[theorem]{Lemma}
\newtheorem{prop}[theorem]{Proposition}
\newtheorem{corol}[theorem]{Corollary}

\allowdisplaybreaks[4]

\begin{document}

\title[Robin $p$\,-Laplacian in~$C^1$~domains]{An eigenvalue estimate for a Robin $p$\,-Laplacian in~$C^1$~domains}

\author{Konstantin Pankrashkin}


\address{Carl von Ossietzky Universit\"at, Fakult\"at V -- Mathematik und Naturwissenschaften,
Institut f\"ur Mathematik, 26111 Oldenburg, Germany}

\email{konstantin.pankrashkin@uol.de}

\begin{abstract}
Let $\Omega\subset \mathbb{R}^n$ be a bounded $C^1$ domain and $p>1$. For $\alpha>0$, define the quantity
\[
\Lambda(\alpha)=\inf_{u\in W^{1,p}(\Omega),\, u\not\equiv 0} \Big(\int_\Omega |\nabla u|^p\,\mathrm{d}x - \alpha \int_{\partial\Omega} |u|^p \,\mathrm{d} s\Big)\Big/ \int_\Omega |u|^p\,\mathrm{d} x
\]
with $\dd s$ being the hypersurface measure, which is the lowest eigenvalue of the $p$-laplacian in $\Omega$
with a non-linear $\alpha$-dependent Robin boundary condition. We show the asymptotics $\Lambda(\alpha)
=(1-p)\alpha^{p/(p-1)}+o(\alpha^{p/(p-1)})$ as $\alpha$ tends to $+\infty$.
The result was only known for the linear case $p=2$ or under stronger smoothness assumptions.
Our proof is much shorter and is based on completely different and elementary arguments,
and it allows for an improved remainder estimate for $C^{1,\lambda}$ domains.
\end{abstract}

\subjclass[2010]{Primary: 35J92, 35P15, 49R05, 49J40, 35J05}
\keywords{Eigenvalue, $p$-Laplacian, Robin boundary condition}


\maketitle


\section{Introduction}

Let $\Omega\subset \RR^n$ be a bounded $C^1$ domain and $p>1$. For $\alpha>0$ define
\[
\Lambda(\alpha):=\inf_{u\in W^{1,p}(\Omega),\, u\not\equiv 0} \Big(\int_\Omega |\nabla u|^p\,\mathrm{d}x - \alpha \int_{\partial\Omega} |u|^p \,\mathrm{d} s\Big)\Big/
\int_\Omega |u|^p\,\mathrm{d} x,
\]
where $\dd s$ is the hypersurface measure. Denote $q:=\frac{p}{p-1}$ the H\"older conjugate of $p$.
In this note we prove the following result:
\begin{theorem}\label{thm1}
As $\alpha$ tends to $+\infty$ there holds
\begin{equation}
  \label{eqla}
\Lambda(\alpha)=(1-p)\alpha^q+o(\alpha^q).
\end{equation}
In addition, if $\Omega$ is of class $C^{1,\lambda}$ with some $\lambda\in(0,1)$, then the remainder estimate can be improved to
$O(\alpha^k)$ with $k:=\frac{q+\lambda}{1+\lambda} \in \big(\frac{q+1}{2}, q\big)$.
\end{theorem}
A standard argument shows that $\Lambda(\alpha)$ can be viewed
as the lowest eigenvalue of a non-linear problem involving the $p$-laplacian
$\Delta_p: u\mapsto \nabla \cdot \big( |\nabla u|^{p-2}\nabla u \big)$,
\[
-\Delta_p u= \Lambda\,  |u|^{p-2}u \text{ in } \Omega, \quad
|\nabla u|^{p-2}\tfrac{\partial u}{\partial \nu}=\alpha |u|^{p-2} u \text{ on } \partial\Omega,
\]
with $\nu$ being the outer unit normal, see e.g. \cite{kp1}. It seems that the asymptotic behavior for large $\alpha$
was first addressed by Lacey, Ockedon and Sabina~\cite{los} for the linear situation ($p=2$), and $C^1$ domains represent a borderline case. On one hand, the asymptotic behavior \eqref{eqla} is not valid for non-smooth domains \cite{bp,khalile2,lp}. On the other hand, for $C^{1,1}$ domains one has
$\Lambda(\alpha)= (1-p)\alpha^q-H \alpha +o(\alpha)$
with $H$ being the maximum mean curvature of the boundary~\cite{kp1},
which is much more detailed, but the proof depends heavily on the existence of a tubular neighborhood of the boundary and on
the regularity of boundary curvatures. We refer to  the review \cite[Sec.~4.4.2]{bfk} for a detailed discussion
of the linear case. Lou and Zhu in~\cite{lz} established the asymptotics \eqref{eqla}
for $C^1$ domains and $p=2$ using an involved
combination of a blow up argument with the non-existence of positive solutions for some linear boundary value problems.
We provide a very short elementary proof based on an integration by parts combined with a simple regularization,
which is by introducing additional ingredients in the proof of the Sobolev trace theorem in Grisvard's book \cite[Sec.~1.5]{gris}.
This allows one to include all values of $p$ and to obtain an improved remainder estimate for $C^{1,\lambda}$
domains with minimal additional effort. Remark that our proof works for $\lambda=1$ as well, but the final result is weaker than the one of \cite{kp1} mentioned above.
By analogy, we do not expect our remainder estimate for $C^{1,\lambda}$ case to be optimal: 
we collect some observations in Section~\ref{remest}.
In particular, Corollary~\ref{corol4} states that the remainder estimate
for  $C^1$ domains cannot be improved to $O(\alpha^r)$ with $r<q$.
 
\section{Proof of Theorem~\ref{thm1}}

\subsection{Lower bound}
We will need a very simple extension result given by the following Lemma~\ref{lem2}. While it is just a version
of Tietze extension theorem, see e.g. \cite[\S8.4]{schub}, we prefer to give a complete direct proof at the end of this subsection.
\begin{lemma}\label{lem2}
The outer unit normal $\nu$ on $\partial\Omega$ can be extended to a continuous map $\mu:\RR^n\to \RR^n$ satisfying $\big|\mu(x)\big|\le 1$ for all $x\in\overline{\Omega}$. If $\Omega$ is of class $C^{1,\lambda}$,
then the above map $\mu$ may be chosen of class $C^{0,\lambda}$.
\end{lemma}

Let $\mu$ as in Lemma~\ref{lem2} and pick $\varepsilon \in(0,1)$, then there exists a $C^1$ map $\mu_\varepsilon:\RR^n\to\RR^n$
with $\big| \mu_\varepsilon(x)-\mu(x)\big|\le \varepsilon$ for all $x\in \overline{\Omega}$.
For any $u\in C^1(\overline{\Omega})$ the divergence theorem gives
\begin{equation}
 \label{2}
\int_{\partial\Omega} |u|^p \mu_\varepsilon\cdot \nu\,\dd s=\int_\Omega \nabla\cdot \big(|u|^p\mu_\varepsilon\big)\dd x=\int_\Omega \big( p \,|u|^{p-2}u\nabla u\cdot \mu_\varepsilon + |u|^p (\nabla\cdot \mu_\varepsilon)\big)\dd x.
\end{equation}
In  $\Omega$ we estimate pointwise $|\mu_\varepsilon|\le |\mu|+\varepsilon\le 1+\varepsilon$ and $|\nabla\cdot \mu_\varepsilon|\le c_\varepsilon$
with a suitable constant $c_\varepsilon>0$. On $\partial\Omega$ we have $\mu=\nu$ and $\mu_\varepsilon\cdot \nu= \nu\cdot \nu +(\mu_\varepsilon-\mu)\cdot \nu\ge 1-\varepsilon$. Therefore, Eq. \eqref{2} yields
\[
(1-\varepsilon)\int_{\partial\Omega} |u|^p \dd s\le (1+\varepsilon)p\int_\Omega |u|^{p-1}|\nabla u|\dd x
+ c_\varepsilon \int_\Omega |u|^p\dd x.
\]
Using Young's inequality, for any $\delta>0$ we estimate
\[
|u|^{p-1}|\nabla u|= \big(\delta^{-1}|u|^{p-1}\big) \big(\delta|\nabla u|\big)
\le \tfrac{1}{q}\,\delta^{-q}|u|^p  + \tfrac{1}{p} \,\delta^p |\nabla u|^p,
\]
then
\[
(1-\varepsilon)\int_{\partial\Omega} |u|^p \dd s\le(1+\varepsilon)\delta^p\int_\Omega |\nabla u|^p \dd x
+ \Big(c_\varepsilon+(1+\varepsilon)(p-1)\delta^{-q}\Big) \int_\Omega |u|^p\dd x,
\]
which holds by density for all $u\in W^{1,p}(\Omega)$ and can be transformed into
\[
\int_\Omega |\nabla u|^p \dd x- \tfrac{1-\varepsilon}{1+\varepsilon} \,\delta^{-p}\int_{\partial\Omega} |u|^p \dd s\ge
-\delta^{-p}\Big( \tfrac{1}{1+\varepsilon}\,c_\varepsilon+(p-1)\delta^{-q}\Big) \int_\Omega |u|^p\dd x.
\]
For $\delta=\big(\tfrac{1+\varepsilon}{1-\varepsilon}\, \alpha\big)^{-1/p}$
it takes the form
\begin{align*}
\int_\Omega |\nabla u|^p \dd x-\alpha \int_{\partial\Omega} |u|^p \dd s &\ge 
\Big((1-p) \big(\tfrac{1+\varepsilon}{1-\varepsilon}\big)^{q} \,\alpha^{q}
-\tfrac{1}{1-\varepsilon}\,c_\varepsilon\alpha\Big) \int_\Omega |u|^p\dd x.
\end{align*}
Therefore,
\begin{gather}
  \label{eqlow1}
\Lambda(\alpha)\ge (1-p) \big(\tfrac{1+\varepsilon}{1-\varepsilon}\big)^{q} \,\alpha^{q}
-\tfrac{1}{1-\varepsilon}\,c_\varepsilon\alpha,
\end{gather}
and then $\liminf_{\alpha\to +\infty } \alpha^{-q} \Lambda(\alpha)  \ge (1-p)\big(\tfrac{1+\varepsilon}{1-\varepsilon}\big)^{q}$.
As $\varepsilon\in(0,1)$ can be taken arbitrarily small, we arrive at $\liminf_{\alpha\to +\infty } \alpha^{-q} \Lambda(\alpha) \ge 1-p$
giving the lower bound.

Now assume that $\Omega$ is of class $C^{1,\lambda}$. By Lemma~\ref{lem2}, the
map $\mu$ in the preceding computations can be assumed $C^{0,\lambda}$. The idea of an improved
remainder estimate is to apply the standard mollifying procedure to construct $\mu_\varepsilon$
and then to control the constant $c_\varepsilon$ in the above computations using the modulus of continuity of $\mu$.
Namely, let $\rho\in C^\infty_c(\RR^n)$ be non-negative, supported in the unit ball centered at the origin, with $\int \rho=1$.
For $t>0$
consider the function $\rho_t:x\mapsto t^{-n}\rho (t^{-1}x)$ and then the $C^\infty $vector field $m_t:=\mu\star \rho_t\in C^\infty(\RR^n,\RR^n)$, with $\star$ being the convolution product. One has
\[
\nabla\cdot m_t(x)=\tfrac{1}{t^{n+1}}\int_{\RR^n} \mu(y) \cdot (\nabla\rho)\big( \tfrac{x-y}{t}\big)\dd y
=\tfrac{1}{t}\int_{|z|<1} \mu(x-tz) \cdot \nabla \rho(z)\dd z,
\]
so that $\|\nabla \cdot m_t\|_{L^\infty(\Omega)}\le a t^{-1}$
with $a:=\|\mu\|_{L^\infty(\Omega_1,\RR^n)}\|\nabla \rho\|_{L^1(\RR^n,\RR^n)}$, where
we denote $\Omega_1:=\{x+z:\, x\in\Omega, \, |z|<1\}$. Furthermore,
\[
m_t(x)-\mu(x)=\int_{\RR^n} \big( \mu(x-y)-\mu(x) \big) \rho_t(y)\dd y=\int_{|z|<1} \big(\mu(x-tz) -\mu (x)\big)\rho(z)\dd z.
\]
As $\mu$ is $C^{0,\lambda}$, with a suitable $b>0$ one 
has $\big|\mu(x-tz) -\mu (x)\big|\le b |tz|^\lambda$ and then $\big| m_t(x)-\mu(x)\big|\le b t^\lambda$
for all $x\in \overline{\Omega}$ and $t>0$. 
Hence, if for $t\in(0,b^{-1/\lambda})$ one sets $\varepsilon:=bt^\lambda\in(0,1)$ and $\mu_\varepsilon:=m_t$,
then in $\overline{\Omega}$ one has $|\mu_\varepsilon-\mu|\le\varepsilon$ and
$|\nabla\cdot \mu_\varepsilon|\le c_\varepsilon:= c \varepsilon^{-1/\lambda}$ with $c:=a b^{1/\lambda}$,
and the inequality~\eqref{eqlow1} takes the form
\[
\Lambda(\alpha)\ge (1-p) \big(\tfrac{1+\varepsilon}{1-\varepsilon}\big)^{q} \,\alpha^{q}
-\tfrac{1}{1-\varepsilon}\,c \varepsilon^{-1/\lambda}\alpha
\text{ for all $\alpha>0$ and $\varepsilon\in(0,1)$.}
\]
Now we take $\varepsilon=\varepsilon(\alpha)$ with $\lim_{\alpha\to+\infty}\varepsilon(\alpha)=0$, then by applying
Taylor expansions one arrives at $\Lambda(\alpha)\ge (1-p)\,\alpha^{q}+O( \varepsilon \alpha^{q}+\varepsilon^{-1/\lambda}\,\alpha)$.
To optimize the last summand we now set  $\varepsilon:=\alpha^{\lambda (1-q)/(\lambda+1)}$, which gives $\Lambda(\alpha)\ge (1-p)\alpha^q+O(\alpha^k)$ with $k:=\frac{q+\lambda}{1+\lambda}$ for large $\alpha$.

\begin{proof}[Proof of Lemma~\ref{lem2}]
The idea is very standard: one first construct an explicit extension near each point of the boundary,
then these constructions are glued together using a partition of unity. The local construction
is also very simple: informally, in a compact coordinate patch over which the boundary
of the domain is given by a graph $y=\phi(x)$,
the local extension can be defined by $\mu(x,y)=\nu\big(x,\phi(x)\big)$.

Let us describe the above procedure in a detailed rigorous way.
By the usual definition of a $C^1$ domain, see e.g. \cite[Def.~1.2.1.1]{gris}, each point of $\partial\Omega$
admits an open neighborhood $V$ with the following properties:
\begin{itemize}
\item  the set $V$ is a hyperparallelepiped, i.e. there exist orthogonal coordinates $y=(y_1,\dots,y_n)$ and strictly positive numbers $a_1,\dots,a_n$ such that
$V=\big\{y: y'\in V',\, y_n\in(-a_n,a_n)\big\}$
with  $V':=(-a_1,a_1)\times\dots\times(-a_{n-1},a_{n-1})$ and $y':=(y_1,\dots,y_{n-1})$,
\item there exists a $C^1$ function $\varphi:V'\to \RR$ with
 $|\varphi(y')|\le \frac{1}{2} \,a_n$ for all $y'\in V'$ and such that
$\Omega\mathbin{\cap} V=\big\{y \in V: y_n<\varphi(y')\big\}$ and $\partial\Omega\mathbin{\cap} V=\big\{y \in V: y_n=\varphi(y')\big\}$.
\end{itemize}
The map $\nu_0: V\ni y\mapsto \big(1+|\nabla\varphi(y')|^2\,\big)^{-\frac{1}{2}}\big(-\partial_1 \varphi(y'),\dots,-\partial_{n-1} \varphi(y'),1\big)\in\RR^n$
is clearly continuous, and for $y\in\partial\Omega\mathbin{\cap} V$, i.e. for $y_n=\varphi(y')$, it coincides with the outer unit normal $\nu$.
Hence, the map $\nu_0$ is a continuous extension of $\nu:\partial\Omega\mathbin{\cap} V\to\RR^n$ to the whole of $V$,
and it satisfies pointwise $|\nu_0|=1$ by construction. If $\Omega$  is of class $C^{1,\lambda}$,
then one can additionally assume that $\varphi\in C^{1,\lambda}(V',\RR)$, and then $\partial_j\varphi\in C^{0,\lambda}(V',\RR)$.
Simple manipulations with H\"older continuous functions, see e.g. \cite[Sec.~1.2]{fior}, show that $\nu_0\in C^{0,\lambda}(V,\RR^n)$.

The boundary $\partial\Omega$ is compact and can be covered by finitely many open hyperparallelepipeds $V_1,\dots,V_m$ with the above properties,
and we denote by $\nu_j$ continuous maps $V_j\to \RR^n$ extending $\nu:\partial\Omega\mathbin{\cap} V_j\to\RR^n$ and satisfying $|\nu_j|\le 1$,
which exist due to the preceding construction. As $\Omega$, $V_1,\dots,V_m$ form a finite open covering
of $\overline{\Omega}$,
one can find a subordinated partition of unity, i.e. $C^\infty$ functions $\psi_0,\dots,\psi_m:\RR^n\to[0,1]$
with $\supp\psi_0\subset \Omega$, $\supp \psi_j\subset V_j$ for $j=1,\dots,m$, and $\psi_0+\dots+\psi_m=1$ in  
$\overline{\Omega}$. We define $\mu:\RR^n\to\RR^n$ by
$\mu:=\psi_1 \nu_1+\dots+\psi_m\nu_m$, which is continuous as each summand is a continuous function.
For $x\in\partial\Omega$ we have $\psi_0(x)=0$ and $(\psi_j\nu_j)(x)=\psi_j(x)\nu(x)$ for $j=1,\dots,m$, which yields
\[
\mu(x)=(\psi_1 \nu_1) (x)+\dots+(\psi_m\nu_m)(x)=\big(\psi_0(x)+\psi_1(x)+\dots\psi_m(x)\big)\,\nu(x)=\nu(x)
\]
and shows that $\mu$ is an extension of $\nu$. Finally, for any $x\in\overline{\Omega}$ we have
\begin{align*}
\big|\mu(x)\big|&\le \big|(\psi_1 \nu_1) (x)\big|+\dots+\big|(\psi_m\nu_m)(x)\big|\\
&=\psi_1(x)\big|\nu_1(x)\big|+\dots+\psi_m(x)\big|\nu_m(x)\big|\\
&\le \psi_1(x)+\dots+\psi_m(x)=1-\psi_0(x)\le 1.
\end{align*}
If $\Omega$ is $C^{1,\lambda}$, then one can assume $\nu_j\in C^{0,\lambda}(V_j,\RR^n)$,
and then $\mu\in C^{0,\lambda}(\RR^n,\RR^n)$.
\end{proof}

\subsection{Upper bound}\label{ssupp}

The upper bound was already obtained in \cite[Prop.~6.2]{kp1}
using a minor variation of a construction by Giorgi and Smits in~\cite[Thm.~2.3]{gs1}.
As the argument is very simple, we repeat it here in order to have a self-contained presentation.
Set $\beta:=\frac{q}{p}$ and consider the function $u:x\mapsto e^{\beta x_1}$ and the vector field
 $F:x\mapsto e^{p\beta x_1} e_1$ with $e_1:=(1,0,\dots,0)$, then
\[
\int_{\partial \Omega} |u|^p\dd s
=\int_{\partial\Omega}e^{p\beta x_1}\dd s
=\int_{\partial\Omega} F\cdot e_1\dd s =
\int_{\partial\Omega} F\cdot \nu \,\dd s + \int_{\partial\Omega} e^{p\beta x_1} (1- e_1\cdot \nu) \,\dd s.
\]
Using the quantity
\[
I(\alpha):=\int_{\partial\Omega} e^{p\beta x_1} (1- e_1\cdot \nu) \,\dd s> 0
\]
and the divergence theorem we arrive at
\[
\int_{\partial \Omega} |u|^p\dd s=\int_{\Omega} \nabla\cdot F\,\dd s +I(\alpha)=p\beta \int_\Omega |u|^p \dd x+I(\alpha),
\]
and then
\begin{equation}
   \label{lka1}
\begin{aligned}
\Lambda(\alpha) &\le \Big(\int_\Omega |\nabla u|^p\,\mathrm{d}x - \alpha \int_{\partial\Omega} |u|^p \,\mathrm{d} s\Big)\Big/
\int_\Omega |u|^p\,\mathrm{d} x\\
&\le 
 \Big(\beta^p \int_\Omega  |u|^p\dd x -p\beta \alpha \int_\Omega |u|^p \dd x - \alpha I(\alpha)\Big) \Big/ \int_\Omega  |u|^p\dd x\\
&=\beta^p -p\beta\alpha-K(\alpha)\equiv (1-p)\alpha^q-K(\alpha), 
\end{aligned}
\end{equation}
where
\begin{equation}
   \label{aka}
K(\alpha):=\alpha\int_{\partial\Omega} e^{p\beta x_1} (1- e_1\cdot \nu) \,\dd s \Big/  \int_\Omega  e^{p\beta x_1}\,\dd x> 0.
\end{equation}
which implies $\Lambda(\alpha)< (1-p)\alpha^q$ and concludes the proof of Theorem~\ref{thm1}.

\section{More on remainder estimates}\label{remest}

We continue using the notation of subsection~\ref{ssupp} and remark that the term 
$K(\alpha)$ was kept intentionally, as this allows one to discuss the optimality of the remainders.
We restrict our attention to two-dimensional domains.
\begin{prop}\label{prop3}
For any $\lambda\in(0,1)$ one can find a bounded $C^{1,\lambda}$ domain $\Omega\subset \RR^2$ and a constant $c_\lambda>0$ such that
for large $\alpha$ there holds
\[
\big| \Lambda(\alpha) - (1-p) \alpha^q\big|\ge c_\lambda \alpha^m, \quad m:=\dfrac{p-\frac{2\lambda}{\lambda+1}}{p-1}\equiv \dfrac{q+(2-q)\lambda}{1+\lambda}\in (1,q).
\]
\end{prop}
\begin{proof}
Let $\delta>0$. As the function $t\mapsto |t|^{1+\lambda}$ is $C^{1,\lambda}$ on any finite interval,
one can find a bounded $C^{1,\lambda}$ domain
$\Omega\subset \RR^2$ such that
\[
\Omega\cap \big\{(x_1,x_2):\, x_1>-\delta\big\}=\big\{ (x_1,x_2): -\delta<x_1<-|x_2|^{1+\lambda}\big\}.
\]
First, there holds
\[
\int_{\partial\Omega} e^{p\beta x_1} (1- e_1\cdot \nu) \,\dd s\ge \int_{\partial\Omega \cap\{x_1> -\delta\}} e^{p\beta x_1} (1- e_1\cdot \nu) \,\dd s.
\]
Using the parametrization $(0,\delta)\ni t\mapsto (-t, \pm t^{1/(1+\lambda)})$ of $ \partial\Omega \cap\{x_1>0, \pm x_2>0\}$
and assuming that $\delta$ is chosen small enough and that $\alpha$ is large enough, we get
\begin{multline*}
\int_{\partial\Omega \cap\{x_1> -\delta\}} e^{p\beta x_1} (1- e_1\cdot \nu) \,\dd s\\
\begin{aligned}
&=2\int_0^\delta e^{-p\beta t} \bigg( 1- \dfrac{t^{-\lambda/(1+\lambda)}}{\sqrt{1+ \big(\tfrac{1}{1+\lambda}\big)^2\,t^{-2\lambda/(1+\lambda)}}}\bigg) \, \sqrt{1+ \big(\tfrac{1}{1+\lambda}\big)^2\,t^{-2\lambda/(1+\lambda)}}\,\dd t\\
&=2\int_0^\delta e^{-p\beta t} \Big( \sqrt{1+ \big(\tfrac{1}{1+\lambda}\big)^2\,t^{-2\lambda/(1+\lambda)}} - t^{-\lambda/(1+\lambda)}\Big) \,\dd t\\
&=\tfrac{2}{1+\lambda}\int_0^\delta e^{-p\beta t} t^{-\lambda/(1+\lambda)} \Big(\sqrt{1+(1+\lambda)^2t^{2\lambda/(1+\lambda)}} - 1\Big)\,\dd t\\
&\ge  \tfrac{2}{1+\lambda}\int_0^\delta e^{-p\beta t} t^{-\lambda/(1+\lambda)} \Big(\dfrac{1}{4}\,(1+\lambda)^2t^{2\lambda/(1+\lambda)}\Big)\, \dd t\\
&= \tfrac{1+\lambda}{2} \int_0^\delta e^{-p\beta t} t^{\lambda/(1+\lambda)} \dd t\ge \tfrac{1+\lambda}{2} \int_0^\infty e^{-p\beta t} t^{\lambda/(1+\lambda)} \dd t -e^{-\delta\beta}\\
&= \tfrac{1+\lambda}{2}\, \Gamma\big( \tfrac{1+2\lambda}{1+\lambda}\big) (p\beta)^{-(1+2\lambda)/(1+\lambda)} -e^{-\delta\beta}.
\end{aligned}
\end{multline*}
Therefore,
\[
\int_{\partial\Omega} e^{p\beta x_1} (1- e_1\cdot \nu) \,\dd s\ge 
\tfrac{1+\lambda}{2}\,
\Gamma\big( \tfrac{1+2\lambda}{1+\lambda}\big) (p\beta)^{-(1+2\lambda)/(1+\lambda)} - e^{-\delta\beta}.
\]
In addition, 
\[
\int_\Omega  e^{p\beta x_1}\,\dd x\le \int_{\Omega \cap\{x_1\in(-\delta,0)\}}  e^{p\beta x_1}\,\dd x +|\Omega| e^{-p\delta\beta},
\]
while the first summand on the right-hand side is estimated as
\begin{multline*}
\int_{\Omega \cap\{x_1> -\delta\}}  e^{p\beta x_1}\,\dd x=\int_{-\delta}^0 \int_{-|x_1|^{1/(1+\lambda)}}^{|x_1|^{1/(1+\lambda)}} e^{p\beta x_1}\, \dd x_2\, \dd x_1\\
\begin{aligned}
&=2\int_0^\delta t^{1/(1+\lambda)} e^{-p \beta t}\dd t\le 2\int_0^\infty t^{1/(1+\lambda)} e^{-p \beta t}\dd t\\
&=2 (p\beta)^{-(2+\lambda)/(1+\lambda)} \int_0^\infty t^{1/(1+\lambda)} e^{-t}\dd t=2\Gamma\big( \tfrac{2+\lambda}{1+\lambda}\big) (p\beta)^{-(2+\lambda)/(1+\lambda)}.
\end{aligned}
\end{multline*}
Putting these estimates into the expression \eqref{aka} for $K(\alpha)$ one sees that for a suitable $c_\lambda>0$ we have, as $\alpha$ is large,
\[
K(\alpha)\ge c_\lambda \alpha \beta^{(1-\lambda)/(1+\lambda)} \equiv c_\lambda \alpha^m, \quad m:=\big(p-\tfrac{2\lambda}{\lambda+1}\big)/(p-1),
\]
and the claim follows from the above inequality \eqref{lka1}.
\end{proof}

\begin{corol}\label{corol4}
For any $r<q$ one can find a bounded $C^1$ domain $\Omega\subset \RR^2$ and a constant $c_r>0$
satisfying $\big|\Lambda(\alpha) - (1-p) \alpha^q\big| \ge c_r \,\alpha^r$ for large $\alpha$.
\end{corol}

\begin{proof}
Given a value of $r$, one can find a sufficiently small  $\lambda\in(0,1)$ to have
$(p-\frac{2\lambda}{\lambda+1})/(p-1)\ge r$. For this value of $\lambda$ one applies Proposition~\ref{prop3}.
\end{proof}

%

\end{document}